\documentclass[12pt]{article}
\usepackage{amsfonts,amsmath,amssymb,bbm}
\usepackage{setspace} 
\def\proof{\noindent{\bf Proof:}\hskip10pt}        
\def\QED{\hfill $\Box$}

\font\tenmath=msbm10 scaled 1200
\font\sevenmath=msbm7 scaled 1200
\font\fivemath=msbm5 scaled 1200
\newfam\mathfam \textfont\mathfam=\tenmath
\scriptfont\mathfam=\sevenmath \scriptscriptfont\mathfam=\fivemath

\begin{document}
\def \\ { \cr }
\def\R{\mathbb{R}}
\def \1{1 \mkern -6mu 1} 
\def\N{\mathbb{N}}
\def\E{\mathbb{E}}
\def\P{\mathbb{P}}
\def\Z{\mathbb{Z}}
\def\Q{\mathbb{Q}}
\def\C{\mathbb{C}}
\def\D{\mathbb{D}}
\def\T{\mathbb{T}}
\def \S{{\mathcal P}_{\tt m}}
\def \e{{\rm e}}
\def \i{{\rm i}}
\def \f{{\mathcal F}}
\def \PM{{\mathcal P}_{1}}
\def \PN{{\mathcal P}_{\N}}
\def \Pn{{\mathcal P}_{[n]}}
\def \g{{\mathcal G}}
\def \h{{\mathcal H}}
\def \d{{\tt d}}
\def \k{{\tt k}}
\def \r{{\mathcal R}}
\newtheorem{theorem}{Theorem}
\newtheorem{definition}{Definition}
\newtheorem{proposition}{Proposition}
\newtheorem{lemma}{Lemma}
\newtheorem{corollary}{Corollary}
\centerline{\LARGE \bf Two-parameter Poisson-Dirichlet}
\vskip 2mm
\centerline{\LARGE \bf   measures and reversible exchangeable }
\vskip 2mm
\centerline{\LARGE \bf  fragmentation-coalescence processes}

\vskip 1cm
\centerline{\Large \bf Jean Bertoin}
\vskip 1cm
\noindent
\centerline{\sl Laboratoire de Probabilit\'es, Universit\'e Pierre et Marie Curie}
\centerline{\sl  and DMA, Ecole Normale Sup\'erieure}
\centerline{\sl 175, rue du Chevaleret} 
\centerline{\sl F-75013 Paris, France}
\vskip 15mm

\noindent{\bf Summary. }{\small  } We show that for $0<\alpha<1$ and $\theta>-\alpha$, the Poisson-Dirichlet distribution with parameter $(\alpha, \theta)$ is the unique reversible distribution of a rather natural fragmentation-coalescence process.
This completes earlier results in the literature for certain split and merge transformations 
and the parameter $\alpha =0$.
 \vskip 5mm
\noindent
{\bf A.M.S. Classification.}  60 J 25, 60 J 27
\vskip 3mm
\noindent{\bf e-mail.} {\tt jbe@ccr.jussieu.fr}

\begin{section}{Introduction}
For every $\theta>0$, consider a sequence $\beta_1, \beta_2, \ldots$ of i.i.d. random variables with the beta$(1,\theta)$ distribution, i.e. $\P(\beta_i\in{\rm d}u)
= \theta(1-u)^{\theta-1}{\rm d}u$ with $u\in]0,1[$. The Residual Allocation Model defined by the sequence of random variables
$$\xi_1=\beta_1\ , \ \xi_2=\beta_2(1-\beta_1)\ , \ \ldots\  , \
\xi_{n+1}= \beta_{n+1} \times \prod_{i=1}^n (1-\beta_i)\ , \ \ldots$$
has the Griffiths-Engen-McClosey distribution with parameter $\theta$, GEM$(\theta)$. Rearranging the terms of the sequence $(\xi_n)_{n\in\N}$ in the decreasing order
yields a random variable ${\xi}^{\downarrow}$ with values in the space of (proper) partitions of the unit mass
$$\PM:=\left\{{\bf x}=(x_1,x_2,\ldots): x_1\geq x_2\geq \cdots\geq 0 \hbox{ and } \sum_{n=1}^{\infty}x_n = 1\right\}\,.$$
The ranked sequence ${\xi}^{\downarrow}$ has a Poisson-Dirichlet distribution, which is denoted here by ${\rm PD}(0,\theta)$.  

The Poisson-Dirichlet laws ${\rm PD}(0,\theta)$ form a one-parameter  family of probability measures on $\PM$
that play an important role in a variety of areas. Let us just mention, first for $\theta = 1$,
limit theorems for large combinatorial structures (e.g. decomposition of permutations into cycles, or of integers into prime factors), and then for general $\theta>0$,
 Bayesian statistics, population genetics ... We refer to
the book \cite{ABT} by Arratia, Barbour and Tavar\'e for precise references and much more on this topic. 
The motivation for the present note stems from the fact that 
${\rm PD}(0,\theta)$ arises as invariant distribution for certain
split and merge transformations; see Tsilevich \cite{Tsi}, Gnedin and Kerov \cite{GnKe}, Pitman \cite{PiGEM}, Diaconis {\it et al.}
\cite{DMWZZ}, Mayer-Wolf {\it et al.} \cite{MWZZ}, Schramm \cite{Sch},~...
Specifically, consider the discrete time Markov chain with values in $\PM$
whose transitions can be described as follows. When the chain starts from some 
configuration ${\bf x}=(x_i)_{i\in\N}$,  pick two i.i.d. random integers 
$N$ and $N'$ according to the law $\P(N=i)=\P(N'=i)=x_i$. If $N\neq N'$,
then merge 
 $x_N$ and $x_{N'}$, otherwise split the term $x_N=x_{N'}$ uniformly at random, 
 and in both cases,   leave the other terms of ${\bf x}$ unchanged. Finally,
reorder the resulting terms in the decreasing order to get a partition of the unit mass. The Poisson-Dirichlet distribution ${\rm PD}(0,1)$ is
  the unique invariant (in fact, reversible) probability measure for this
 split and merge evolution. More generally, any ${\rm PD}(0,\theta)$ for $\theta>0$ appears as the invariant law of a variation of the preceding dynamics.
 In this direction, we mention that the question of the existence of reversible laws for 
a large class of coagulation-fragmentation processes with values in certain {\it finite} subsets of $\PM$ (see \cite{Ke, Wh} for some literature on this topic) has been solved
by Durrett {\it et al.} \cite{DGG} and Erlihson and  Granovsky \cite{EG}. However, due to the requirement of finiteness, it does not seem that the latter results can be applied to tackle the problem of existence and uniqueness of the invariant distribution for the split and merge transformation above.

Perman, Pitman and Yor \cite{PPY} extended the notion of Poisson-Dirichlet distributions 
to a  two-parameter family ${\rm PD}(\alpha,\theta)$ with $0\leq \alpha < 1$ and $\theta>-\alpha$. One quick way for this is to apply the Residual Allocation Model above to the situation where $\beta_1, \beta_2, \ldots$ are again independent beta variables, but now $\beta_n$ has parameter $(1-\alpha, \theta+n\alpha)$, viz.
$$\P(\beta_n\in {\rm d}u)={\rm B}(1-\alpha, \theta+n\alpha)^{-1} u^{-\alpha}(1-u)^{\theta+n\alpha-1}{\rm d}u\quad \hbox{ for }0<u<1. $$
 This two-parameter family 
appears naturally in connection with excursion lengths of Bessel processes and their bridges; we refer to the survey by Pitman and Yor  \cite{PPY} 
for a detailed account on this topic.
The purpose of this note is to point out that for any  $0<\alpha < 1$ and $\theta>-\alpha$,
the Poisson-Dirichlet distribution ${\rm PD}(\alpha,\theta)$ 
also arises as the unique reversible law for a certain quite natural fragmentation-coalescence process.

The fragmentation-coalescence processes that we shall deal with here are not directly related to the split and merge transformation, but rather belong to the family of Exchangeable Fragmentation-Coalescence  (EFC) processes which have been introduced recently by Berestycki \cite{Ber}. Roughly speaking, EFC processes evolve in continuous time and combine the dynamics of coalescents with  simultaneous multiple collisions (see \cite{MS, Schw} or Chapter 4 in \cite{RFC})
and that of homogeneous fragmentations (see \cite{Be} or Chapter 3 in \cite{RFC}).
We stress that in general, they are not step Markov processes, in the sense that their first exit from an
arbitrary  configuration is instantaneous a.s.
In this work, collisions in coagulation events are multiple (in fact, every coagulation event involves an infinite number of terms a.s.) and coagulations are simple, in the sense that all the terms involved in a coagulation event merge in a single term. More precisely, using 
the terminology of Pitman \cite{Pimulti}
(see also Sagitov \cite{Sagit}), the rates of multiple collisions are governed by 
the measure $\Lambda_{\theta/\alpha}({\rm d}u)= (1-u)^{\theta/\alpha}{\rm d}u$ for $u\in]0,1[$.
The rates of dislocations are in turn governed by an infinite measure
of Poisson-Dirichlet type denoted by ${\rm PD}(\alpha, -\alpha)$, which has been introduced recently by Basdevant \cite{Bas} in a study of Ruelle's cascades 
 (note that
$(\alpha,-\alpha)$ is excluded from the range of parameters for Poisson-Dirichlet probability measures).
It is interesting to observe that the coagulation rates only depend on the ratio
$\theta/\alpha$ whereas those for dislocations only depend on the parameter $\alpha$.
\end{section}

\begin{section}{Some background}
The purpose of this section is to provide a rigorous presentation of  notions which have been discussed informally above. We thus now introduce the material which will be needed to state and prove our main result in the next section. Exchangeability is one of the major ingredients. We first recall the key connexion between random mass-partitions and random exchangeable partitions of $\N$ due to Kingman, and an important sampling formula  for Poisson-Dirichlet distributions due to Pitman (we refer to e.g. Section 2.3 in \cite{RFC} for a complete account). Then we  briefly  present some main features on EFC processes which are gleaned from \cite{Ber}.

\subsection{Exchangeable random partitions}
A partition $\pi$ of $\N:=\{1,2,\ldots\}$ is an infinite sequence of pairwise disjoint blocks $\pi_1, \pi_2, \ldots$ which are ordered according to the increasing order of their least element, and such that $\bigcup_{i\in\N}\pi_i=\N$. 
The space of partitions of $\N$ is denoted by $\PN$.
 For every $n\in\N$, the restriction to  $[n]:=\{1,\ldots, n\}$ induces
 is a natural projection $\pi\to \pi_{\mid [n]}$ from $\PN$
to  $\Pn$, the space of partitions of $[n]$. Plainly, the sequence of restricted partitions
$(\pi_{\mid [n]}, n\in\N)$ is compatible, in the sense that for every integer $n\geq 2$, $\pi_{\mid [n]}$
is a partition of $[n]$ whose restriction to $[n-1]$ coincides with $\pi_{\mid [n-1]}$. Conversely,
for any sequence of compatible partitions $\pi^{[n]}\in\Pn$ for $n\in\N$, there exists a unique partition $\pi\in\PN$ such that $\pi_{\mid [n]}=\pi^{[n]}$.

A random partition of $\N$ is called {\it exchangeable} when its distribution is invariant by the action of permutations on $\PN$. One can associate to every partition of the unit mass ${\bf x}=(x_1,\ldots)$
an exchangeable random partition $\pi$ by the following construction, which is often referred to as the paint-box process. The partition of the unit mass ${\bf x}$ is viewed as a (proper) probability measure on $\N$. One
introduces a sequence of i.i.d. discrete variables  $\eta_1, \eta_2,\ldots$ distributed according to ${\bf x}$, i.e. $\P(\eta_n=i)=x_i$, and then defines a random partition $\pi$ of $\N$ by declaring that two integers 
$k,\ell$ belong to the same block of $\pi$ if and only if $\eta_k=\eta_\ell$.
Note that $\pi$ contains no singletons a.s., and, by the strong law of large numbers, that every block of $\pi$, say $B$, possesses an asymptotic frequency
$|B|:=\lim_{n\to \infty} n^{-1}{\rm Card}(B\cap [n])$. Further, the sequence $|\pi|=(|\pi_i|,i\in\N)$ 
of the asymptotic frequencies of the blocks of $\pi$ is a size-biased reordering of the sequence
${\bf x}$. Conversely, any exchangeable random partition $\pi$ which contains no singleton a.s. is distributed as a mixture of such  paint-boxes. In particular, all the blocks of $\pi$ possess asymptotic frequencies a.s.,  and their rearrangement in the decreasing order yields a random  (proper) partition
 $|\pi|^{\downarrow}={\bf x}$ of the unit mass whose law characterizes that of $\pi$.

The distribution
 of any exchangeable random partition $\pi$ is determined by an Exchangeable Partition Probability Function (EPPF) 
 $${\tt p}: \bigcup_{k\in\N}\N^k\to [0,1]Ê\quad , \quad {\tt p}(n_1,\ldots, n_k)= \P(\pi_{\mid [n]} = \gamma(n_1,\ldots, n_k))\,,$$
where $n=n_1+\cdots + n_k$ and $\gamma(n_1,\ldots, n_k)$ denotes an arbitrary partition of $[n]$ 
which has $k$ non-empty blocks with cardinals $n_1,\ldots, n_k$ (by the exchangeability property, the preceding probability does not depend on the choice for such a partition $\gamma(n_1,\ldots, n_k)$,
and in particular ${\tt p}$ is a symmetric function of the variables  $n_1,\ldots, n_k$).

We next focus on Poisson-Dirichlet distributions, and fix a pair of parameters $(\alpha, \theta)$ with
 $0<\alpha<1$ and $\theta>-\alpha$. Let ${\bf x}_{\alpha,\theta}$ be a random partition of the unit mass with law ${\rm PD}(\alpha, \theta)$ and $\pi_{\alpha,\theta}$ a random partition of $\N$
 distributed as a paint-box based on ${\bf x}_{\alpha,\theta}$. The EPPF ${\tt p}_{\alpha,\theta}$ of the ${\rm PD}(\alpha,\theta)$ random partition $\pi_{\alpha,\theta}$ is given by Pitman's sampling formula : 
\begin{equation}\label{PSF}
{\tt p}_{\alpha,\theta}(n_1,\ldots,n_k)\,=\,{(\theta/\alpha)_{k\uparrow}\over 
(\theta)_{n\uparrow}}\prod_{i=1}^k-(-\alpha)_{n_i\uparrow}\,,
\end{equation}
where for every
integer $\ell$ and real number $a\not\in \Z_-$,
$$(a)_{\ell\uparrow}=\Gamma(a+\ell)/\Gamma(a)=a(a+1)\cdots (a+\ell-1),$$
and for $\theta=0$, the ratio in the right-hand side of 
\eqref{PSF} is understood as $(k-1)!/(n-1)!$.

\subsection{Exchangeable Fragmentation-Coagulation processes}
This section essentially follows Berestycki \cite{Ber}, with some minor modifications which are made to adapt the presentation to our purposes. 

A Markov process $\Pi=(\Pi(t), t\geq 0)$ with values in $\PN$
is called an EFC process if it is exchangeable (i.e. its law is invariant by the action of permutations on $\PN$-valued processes), and if for every integer $n$, its restriction $\Pi_{\mid [n]}$ to $\Pn$ is a continuous-time Markov chain that only evolves by fragmentation of one block or by coagulation. 
This means that when $\Pi_{\mid [n]}$ starts from some configuration, say $\pi^{[n]}\in\Pn$,
then the only sites at which it can make its first jump are those that can be obtained from $\pi^{[n]}$
by splitting exactly one of its blocks (which thus has at least two elements),  or by merging  at least two of its non-empty blocks. It is important to observe that  an EFC process $\Pi$
can be recovered from the family of its restrictions $\Pi_{\mid [n]}$, and thus its evolution is entirely determined 
by the jump rates of a family of continuous-time Markov chains with values in some finite space, even though the process $\Pi$ itself is not a step process in general.
It is easily seen that, with the obvious exceptions of the cases when the fragmentation or coagulation component is degenerate, the restrictions $\Pi_{\mid [n]}$ are irreducible and thus possess a unique invariant distribution.

The evolution of an EFC process may be easier to understand when the latter is viewed as an infinite particle system, in which particles correspond to blocks of the random partitions $\Pi(t)$. 
Specifically, as $\Pi$ is an exchangeable process,  for every $t\geq 0$, all the blocks of the partition $\Pi(t)$ have asymptotic frequencies 
a.s.\footnote{In fact, Berestycki \cite{Ber} showed the much stronger result that 
with probability one, the partitions $\Pi(t)$ have asymptotic frequencies simultaneously for all $t\geq 0$.}
We write ${\bf X}(t)$ for the partition of the unit mass given by the sequence ranked in the decreasing order of the asymptotic frequencies of the blocks of $\Pi(t)$. We think of the process $({\bf X}(t), t\geq 0)$ as a particle system in which particles are characterized by their masses, and may split or coagulate with other particles.
According to Theorem 7 of Berestycki \cite{Ber},  this particle system has the Feller property, and it is easily seen that  its distribution determines that of $\Pi$.

One of the basic results about EFC processes (cf. Propositions 4 and 5  in \cite{Ber}) states that their dynamics are entirely characterized by a pair of measures on the space of mass-partitions, which specify
respectively the rates of multiple coagulation and of splittings, and 
two nonnegative coefficients, one for the rate of binary coagulations and the other for the erosion. Here, we shall
only consider the situation where there are no binary coagulation and no erosion,
so the aforementioned coefficients are both zero. Further, we shall focus on the case when coagulation events are always simple, i.e. at any given time, all the particles involved into a coagulation event merge in a single new particle. This enables one to represent the measure characterizing  coagulation rates as  some finite measure $\Lambda$ on $]0,1]$. The measure on $\PM$ that characterizes splitting rates will be denoted by $\nu$, it fulfills the conditions
\begin{equation}\label{conddis}
\nu(\{(1,0,\ldots)\})=0\ \hbox{and}\ 
\int_{\PM}(1-x_1)\nu({\rm d}{\bf x})<\infty .
\end{equation}
Informally, the dynamics of the particle system ${\bf X}$ (and hence also of $\Pi$) are determined by the measures
$\Lambda$ and $\nu$ as follows. Coagulation events involving a proportion $u\in]0,1]$ of particles (that are chosen uniformly at random)
occur with rate $u^{-2}\Lambda({\rm d}u)$, and each particle, say $y$, splits into a sequence of smaller particles $yx_1, yx_2, \ldots$ with rate $\nu({\rm d}{\bf x})$, independently of the other particles in the system.
  Conversely, given an arbitrary finite measure $\Lambda$ on $]0,1]$ and an arbitrary measure $\nu$ on $\PM$ such that \eqref{conddis} holds, there exists an EFC process having only simple (but not binary) coagulations, with characteristics $\Lambda$ and $\nu$.

This description of the evolution of the particle system ${\bf X}$ may be
intuitively  appealing, however it would be difficult to make it formal because coagulation (respectively, splitting) events occurs immediately whenever $\int u^{-2}\Lambda({\rm d}u)=\infty $ (respectively, $\nu(\PM)=\infty)$.
 Working with the $\PN$-valued process $\Pi$ and its restrictions
enables one to circumvent this difficulty.
Indeed, the transitions of the continuous time Markov chains $\Pi_{\mid [n]}$
are defined rigorously in terms of $\Lambda$ and $\nu$ as follows.

We first consider coagulation events. 
 For all integers $k,\ell, n$ with $2\leq k \leq \ell\leq n$, pick an arbitrary partition $\gamma\in\Pn$ which has $\ell$ non-empty blocks, and $\gamma'\in\Pn$ another partition which is obtained from
$\gamma$ by a simple coagulation event involving $k$ of its blocks. Then the jump rate of  $\Pi_{\mid [n]}$ corresponding to the coagulation event which transforms
$\gamma$ into $\gamma'$ is
$${\tt c}_n(\gamma,\gamma')={\tt c}(\ell,k)=\int_{[0,1]}u^{k-2}(1-u)^{\ell-k}\Lambda({\rm d}u)\,.$$
These are the only non-zero rates for coagulation events.

We next turn our attention to splitting events. Pick an arbitrary partition $\gamma\in\Pn$ which is not the partition into singletons, and select an arbitrary block $B=(i_1, \ldots, i_k)$ of $\gamma$ with $k\geq 2$ elements.
Then pick an arbitrary partition $\eta\in{\mathcal P}_{[k]}$ which is not trivial (i.e. $\eta\neq ([k], \varnothing, \ldots))$, and split $B$ according to $\eta$.
This means that we replace  the block  $B$ by smaller blocks of the type
$\{i_j, j\in C\}$ where $C$ stands for a generic block of $\eta$, and leave the other blocks
of $\gamma$ unchanged. We write $\gamma'$ for the resulting partition of $[n]$.
Then the jump rate of 
 $\Pi_{\mid [n]}$ corresponding to the splitting event which transforms
 $\gamma$ into $\gamma'$ is
\begin{equation}\label{FR}
{\tt s}_n(\gamma,\gamma')={\tt s}(\eta)=\int_{\PM}\nu({\rm d}{\bf x}) \P(\pi^{\bf x}_{\mid [k]}=\eta)\,,
\end{equation}
where $\pi^{\bf x}$ stands for an exchangeable random partition of $\N$ constructed as a paint-box based on
${\bf x}$. These are the only non-zero rates for splitting events.
It is often convenient to rewrite \eqref{FR} in terms of EPPF. Specifically, for every
partition ${\bf x}$ of a unit mass, write  ${\tt p}_{\bf x}$ for the EPPF of a paint-box based on ${\bf x}$.
Then
\begin{equation}\label{FRbis}
{\tt s}(\eta)={\tt s}(k_1, \ldots, k_{\ell})=\int_{\PM}\nu({\rm d}{\bf x}) {\tt p}_{\bf x}(k_1, \ldots, k_{\ell})\,,
\end{equation}
where $k_1, \ldots, k_{\ell}$ are the cardinals of the (non-empty) blocks of $\eta$.

\end{section}

\begin{section} {Some reversible EFC processes}
Here, we fix two parameters $0<\alpha<1$ and $\theta>-\alpha$ and consider first the finite measure on $]0,1]$ given by
$$\Lambda_{\theta/\alpha}({\rm d}u)=(1-u)^{\theta/\alpha}{\rm d}u.$$
Next, following Basdevant \cite{Bas}, we introduce an infinite measure
${\rm PD}(\alpha,-\alpha)$ on $\PM$ as follows. We denote as usual the generic partition of the unit mass by
${\bf x}=(x_1, x_2, \ldots)$,
and recall that the limit
\begin{equation}
\label{defL}
L=L({\bf x})=\lim_{n\to\infty} n x_n^{\alpha}
\end{equation}
exists ${\rm PD}(\alpha, 0)$-a.s. Further, it is well-known that for every real number $a$,
$$\E_{\alpha,0}(L^a)<\infty \ \Longleftrightarrow \ a>-\alpha\,,$$
where the notation $\E_{\alpha,0}$ indicates that the expectation is computed under
the law ${\rm PD}(\alpha, 0)$. In this setting, the Poisson-Dirichlet
distribution ${\rm PD}(\alpha, \theta)$ is absolutely continuous with respect to
${\rm PD}(\alpha, 0)$ :
$${\rm PD}(\alpha,\theta)=\frac{L^{\theta/\alpha}}
{\E_{\alpha,0}(L^{\theta/\alpha})}{\rm PD}(\alpha,0)\,.$$
Even though $\E_{\alpha,0}(1/L)=\infty$, one can introduce the infinite measure
\begin{equation}\label{PD(a,-a)}
{\rm PD}(\alpha,-\alpha)=L^{-1}{\rm PD}(\alpha,0)\,.
\end{equation}
More precisely, Basdevant \cite{Bas} constructed the measure  ${\rm PD}(\alpha,-\alpha)$ by using a sequence of independent beta variables and a variation of the Residual Allocation Model, and showed in Theorems 4.4 and 4.6 of \cite{Bas} that these two constructions are equivalent. Further ${\rm PD}(\alpha,-\alpha)$ fulfills the requirement
\eqref{conddis}.

Let  ${\Pi}_{\alpha,\theta}$ denote an EFC process with characteristic measures
$\Lambda=\Lambda_{\theta/\alpha}$ and $\nu={\rm PD}(\alpha,-\alpha)$, and let
${\bf X}_{\alpha,\theta}$ stand for the particle system given by the process of the asymptotic frequencies of blocks of ${\Pi}_{\alpha,\theta}$ ranked in the decreasing order.
We are now able to state the main result of this note.

\begin{theorem}\label{T1}
 The Poisson-Dirichlet distribution  ${\rm PD}(\alpha,\theta)$
 is the unique reversible probability measure for the process ${\bf X}_{\alpha,\theta}$.
\end{theorem}

\proof   Fix an integer $n\geq 2$, and denote by $\rho_{\alpha,\theta}^{[n]}$ the probability measure on $\Pn$ which is given by the image of ${\rm PD}(\alpha,\theta)$ by the paint-box process restricted to $[n]$. Pick positive integers $k<n$ and $n_1,\ldots, n_k$
such that $n_1+\cdots + n_k=n$, and consider an arbitrary partition $\gamma\in\Pn$ which has $k$ non-empty blocks with sizes $n_1,\ldots, n_k$.
According to Pitman's sampling formula \eqref{PSF}, we have
$$\rho_{\alpha,\theta}^{[n]}(\{\gamma\})\,=\,{(\theta/\alpha)_{k\uparrow}\over 
(\theta)_{n\uparrow}}\prod_{j=1}^k-(-\alpha)_{n_j\uparrow}\,.$$

Then pick some integer $i\leq k$ such that $n_i\geq 2$. Consider some partition
$\eta$ of $[n_i]$ with $\ell\in\{2,\ldots, n_i\}$ non-empty blocks, say with sizes
$n_{i1},\ldots, n_{i\ell}$, and denote by
$\gamma'\in\Pn$ the partition of $[n]$ which is obtained from $\gamma$
by splitting its $i$-th block according to $\eta$. In particular $\gamma'$ has
$k+\ell-1$ non-empty blocks with respective sizes
$n_1,\ldots, n_{i-1}, n_{i1},\ldots, n_{i\ell},n_{i+1}, \ldots, n_k$.
Another application of Pitman's sampling formula \eqref{PSF} gives
$$\rho_{\alpha,\theta}^{[n]}(\{\gamma'\})\,=\,{(\theta/\alpha)_{k+\ell-1\uparrow}\over 
(\theta)_{n\uparrow}}\left(\prod_{j=1}^k-(-\alpha)_{n_j\uparrow}\right)
\left(\prod_{m=1}^\ell-(-\alpha)_{n_{im}\uparrow}\right) (-(-\alpha)_{n_i\uparrow})^{-1}
\,.$$

On the one hand, the sampling formula of Proposition 4.3 in Basdevant \cite{Bas} 
states that for $\nu={\rm PD}(\alpha,-\alpha)$, the splitting rate ${\tt s}_n(\gamma,\gamma')={\tt s}(\eta)$ defined by \eqref{FR} and \eqref{FRbis} is given explicitly by
$${\tt s}_n(\gamma,\gamma')\,=\,{(\ell-2)!\over -
(-\alpha)_{n_i\uparrow}}\prod_{m=1}^\ell-(-\alpha)_{n_{im}\uparrow}\,.$$
On the other hand, there is a unique way of recovering
 $\gamma$ from $\gamma'$ by a coagulation event. The latter involves exactly $\ell$ blocks amongst the $k+\ell-1$ blocks of $\gamma'$, and thus this coagulation event
occurs with rate
\begin{eqnarray*}{\tt c}_n(\gamma',\gamma)={\tt c}(k+\ell-1, \ell)&=&
\int_{[0,1]}u^{\ell-2}(1-u)^{k-1}\Lambda_{\theta/\alpha}({\rm d}u)\\
&=&\int_{[0,1]}u^{\ell-2}(1-u)^{k-1+\theta/\alpha}{\rm d}u\\
&=&\,\frac{\Gamma(\ell-1)\Gamma(k+\theta/\alpha)}{\Gamma(k+\ell-1+\theta/\alpha)}\\
&=&\frac{(\ell-2)! (\theta/\alpha)_{k \uparrow}}{ (\theta/\alpha)_{ k+\ell-1\uparrow}}.
\end{eqnarray*}

Putting the pieces together, we realize that there is the identity
\begin{equation}\label{dbc}
\rho_{\alpha,\theta}^{[n]}(\{\gamma\}){\tt s}_n(\gamma,\gamma')=
 \rho_{\alpha,\theta}^{[n]}(\{\gamma'\}){\tt c}_n(\gamma',\gamma).
\end{equation}
In other words, for every integer $n\geq 2$, the detailed balance equation 
holds for the Markov chain $\Pi_{\mid [n]}$ and the probability measure $\rho_{\alpha,\theta}^{[n]}$.  It should now be plain from the results which have been recalled in Section 2 that this proves our claim. \QED

Our proof of Theorem \ref{T1} amounts to checking the detailed balance equation \eqref{dbc} by  inspection and does not make the statement quite intuitive.
So it may be interesting to explain quickly some hints for this result.
First, for $\theta=0$, Pitman \cite{Pimulti} showed that the uniform measure $\Lambda_0$ on $[0,1]$ governs the coagulation rates of the Bolthausen-Sznitman coalescent $C^{\rm BS}$. Further, the logarithmic time-change $t\to -\ln t$ transforms 
 the Bolthausen-Sznitman coalescent into a time-inhomogeneous fragmentation
 process $F^{\rm R}(t)=C^{\rm BS}(-\ln t)$, $t\in]0,1[$, which is closely related to Ruelle's cascades. As Basdevant \cite{Bas} identified the instantaneous dislocation measure of
 $F^{\rm R}$ at time $t$ as $t^{-1}{\rm PD}(t,-t)$, this suggests that combining fragmentation with rates governed by ${\rm PD}(t,-t)$ and coagulation with rates
 governed by $\Lambda_0$ might yield a stationary Markov process whose invariant law should be given by the distribution of  the Bolthausen-Sznitman coalescent at time $-\ln t$, that is 
 ${\rm PD}(t,0)$.
For a general $\theta>-\alpha$,  Theorem \ref{T1} has its roots in the remakable duality identity between coagulation and fragmentation operators based on Poisson-Dirichlet variables due to Pitman (see Theorem 12 in \cite{Pimulti} or Theorem 4.4 in \cite{RFC}),  combined with intuition gained from the preceding analysis in the case $\theta=0$.

We now complete this work by proving the convergence to equilibrium for the EFC process ${\bf X}_{\alpha,\theta}$.

\begin{corollary} \label{C1}
 For every configuration ${\bf x}\in \PM$, the distribution of 
${\bf X}_{\alpha,\theta}(t)$ given ${\bf X}_{\alpha,\theta}(0)={\bf x}$
converges weakly to ${\rm PD}(\alpha, \theta)$ when $t\to\infty$.
\end{corollary}

\proof We consider the process $\Pi_{\alpha,\theta}$ with values in $\PN$ started from
an exchangeable random partition given by a paint-box process based on 
${\bf x}$. In particular, for every $t\geq 0$,  the law of  $|\Pi_{\alpha,\theta}(t)|^{\downarrow}$, the sequence of the asymptotic frequencies of the blocks of $\Pi_{\alpha,\theta}(t)$ ranked in the decreasing order, is that
of ${\bf X}_{\alpha,\theta}(t)$ given ${\bf X}_{\alpha,\theta}(0)={\bf x}$.
According to Theorem 8 of Berestycki \cite{Ber}, we know that
$\Pi(t)$ converges in distribution as $t\to \infty$ to the equilibrium measure
of $\Pi$, say $\rho_{\alpha,\theta}$. An application of Proposition 2.9 in \cite{RFC} shows that given ${\bf X}_{\alpha,\theta}(0)={\bf x}$, ${\bf X}_{\alpha,\theta}(t)$  converges in distribution as $t\to \infty$ to the image of $\rho_{\alpha,\theta}$ by the map
$\pi\to |\pi|^{\downarrow}$, which is the stationary law for the EFC process ${\bf X}_{\alpha,\theta}$. By Theorem \ref{T1}, we know that the latter is ${\rm PD}(\alpha,\theta)$.  \QED

We stress that, more generally, this argument can be applied to general EFC processes
to establish the convergence to equilibrium for the mass-partition valued process ${\bf X}$.

\end{section}


\end{document}